\documentclass{amsart}

\usepackage{amssymb}
\usepackage{graphicx}
\usepackage{MnSymbol}

\usepackage{amsthm}

\title[Correction: Root extraction]{Correction:  Root extraction in one-relator groups and slenderness}
\author{Samuel M. Corson}

\theoremstyle{definition}
\theoremstyle{definition}

\theoremstyle{definition}

\theoremstyle{definition}

\theoremstyle{definition}
\theoremstyle{definition}
\theoremstyle{definition}
\theoremstyle{definition}\newtheorem{question}{Question}
\theoremstyle{definition}
\theoremstyle{definition}
\theoremstyle{definition}
\theoremstyle{definition}
\theoremstyle{definition}
\theoremstyle{definition}
\theoremstyle{definition}
\theoremstyle{definition}
\newtheorem*{question*}{Question}
\newtheorem*{theorem*}{Theorem}
\newtheorem*{corollary*}{Corollary}
\newtheorem*{lemma*}{Lemma}
\newtheorem*{fact*}{Fact}
\newtheorem*{definition*}{Definition}

\def\pmc#1{\setbox0=\hbox{#1}
    \kern-.1em\copy0\kern-\wd0
    \kern.1em\copy0\kern-\wd0}

\newcommand{\HEG}{\textbf{HEG}}

\begin{document}

\address{School of Mathematics, University of Bristol, Fry Building, Woodland Road, Bristol, BS8 1UG, United Kingdom.}
\email{sammyc973@gmail.com}
\subjclass[2010]{20F05, 57M07 }
\thanks{This work was supported by the Additional Funding Programme for Mathematical Sciences, delivered by EPSRC (EP/V521917/1) and the Heilbronn Institute for Mathematical Research}

\maketitle

\begin{abstract}  This very short correction notes a gap in an argument of an earlier paper, and also provides a theorem of similar flavor to the main result of that paper.

\end{abstract}

I am indebted to Dawid Kielak \cite{Kie} for pointing out a gap in the proof of Lemma 2.4 in \cite{Co}, which was used to prove Theorems A and B of that paper.  There are no known counterexamples to those results, so they may be regarded as open problems.  The problem with that original argument is that if $L$ and $K$ are two Magnus subgroups of a one-relator group, then the expression of an element in the intersection $L \cap K$ as a word in the free generating set of $L$ can be different from that in the free generating set for $K$.   I'll state and prove a result similar to Theorem B with stronger conclusion but an extra assumption (that the one-relator group includes no Baumslag-Solitar subgroups).  The proof similarly utilizes root extraction.

The notation will be identical with that of the paper \cite{Co}, with $\HEG$, $\HEG^m$, etc., as therein.  We'll recall the following definitions from \cite{ConCor}.

\begin{definition*}  A group $G$ is cm-slender (respectively lcH-slender) if every abstract group homomorphism $\phi: H \rightarrow G$, where $H$ is a completely metrizable (resp. locally compact Hausdorff) topological group, has open kernel.
\end{definition*}

\begin{theorem*}  Let $G$ be a (possibly uncountable) one-relator group which has no Baumslag-Solitar subgroup.  The following hold.

\begin{enumerate}

\item If $\phi: \HEG \rightarrow G$ is an abstract homomorphism then for some $m \in \mathbb{N}$ the image $\phi(\HEG^m)$ is finite.

\item If $\phi: H \rightarrow G$ is an abstract homomorphism, with $H$ either a completely metrizable or locally compact Hausdorff topological group, then there is a normal open subgroup $V \leq H$ with $\phi(V)$ finite.

\end{enumerate}
In particular, a torsion-free one-relator group without Baumslag-Solitar subgroups is n-slender, cm-slender, and lcH-slender.

\end{theorem*}

\begin{proof}  We will first prove claim (1) and then give the quite analogous argument for (2).  Assume that $\phi: \HEG \rightarrow G$ is an abstract group homomorphism, where $G$ is a one-relator group without Baumslag-Solitar subgroups.  If $G = \langle X \mid r\rangle$ we let $Y \subseteq X$ be the set of generators used in the word $r$.  The group $G$ is isomorphic in the natural way to the free product $\langle Y \mid r\rangle * F(X \setminus Y)$ where $F(X \setminus Y)$ is the free group on generators $X \setminus Y$.  Then $\phi: \HEG \rightarrow \langle Y \mid r\rangle * F(X \setminus Y)$, so by \cite[Theorem 1.3]{Ed} we know that for some $N \in \mathbb{N}$ the image $\phi(\HEG^N)$ is included into a conjugate of either $\langle Y \mid r \rangle$ or $F(X \setminus Y)$.  Thus without loss of generality we compose $\phi$ with conjugation by an appropriate element in $G$ so that $\phi(\HEG^N) \leq \langle Y \mid r\rangle$ or $\phi(\HEG^N) \leq F(X \setminus Y)$.  In case $\phi(\HEG^N) \leq F(X \setminus Y)$, since free groups are n-slender \cite[Corollary 3.7]{Edearly}, we can select $m > N$ such that $\HEG^m \leq \ker(\phi \upharpoonright \HEG^N) \leq \ker(\phi)$, so $\phi(\HEG^m)$ is trivial, hence finite.

Suppose now that $\phi(\HEG^N) \leq \langle Y \mid r\rangle$.  In case $J = \langle Y \mid r\rangle$ has torsion we know it is hyperbolic (by the Spelling Theorem of Newman \cite{Ne}).  Then there exists some $m > N$ such that $\phi(\HEG^m)$ is finite \cite[Theorem B]{BC}.  Therefore we may assume that $J$ is torsion-free.  As $G$ does not have Baumslag-Solitar subgroups, we know that $J$ is commutative transitive \cite[Theorem 1.3]{FMRR} and so has unique root extraction (i.e. if $s_0^t = s_1^t$, $t > 0$, then $s_0 = s_1$).  Letting $p$ be a prime greater than the length of the relator $r$, we have that for each nontrivial $s \in J$ there is some $n_s \in \mathbb{N}$ such that the equation $x^{p^{n_s}} = s$ has no solution in $J$ \cite[Theorem 1]{Ne}.  Then by unique root extraction we have that for nontrivial $s \in \langle Y \mid r\rangle$, the set $\{x \in J : (\exists k \in \mathbb{N}) x^{p^k} = s\}$ has cardinality at most $n_s$.  Then in the terminology of \cite{ConCor} the group $J$ has finite $p$-antecedents, and as $J$ is countable and torsion-free, we know that $J$ is n-slender \cite[Theorems A, B(c)]{ConCor}.  Then there exists some $m > N$ such that $\phi \upharpoonright \HEG^m$ is trivial and we have considered the last case for (1).

Now we'll prove (2).  Letting $\phi: H \rightarrow G$ with $H$ completely metrizable (respectively locally compact Hausdorff) and $G = \langle Y \mid r \rangle * F(X \setminus Y)$, with $Y$ finite, we have that either $\ker(\phi)$ is open or $\phi(H)$ lies entirely in a conjugate of $\langle Y \mid r \rangle$ or of $F(X \setminus Y)$ by \cite{Sl} (resp.\cite{MN}).  In case $\ker(\phi)$ is open we are already done.  Else we conjugate appropriately so that without loss of generality either $\phi(H) \leq \langle Y \mid r \rangle$ or $\phi(H) \leq F(X \setminus Y)$.  If $\phi(H) \leq F(X \setminus Y)$ then since free groups are cm-slender \cite{Du} (resp. lcH-slender, also \cite{Du}) we see again that $\ker(\phi)$ is open.  We are left with the case where $\phi(H) \leq \langle Y \mid r \rangle$.  If the group $J = \langle Y \mid r \rangle$ has torsion then it is hyperbolic and by \cite[Theorem A]{BC} there is an open normal subgroup $V \leq H$ such that $\phi(V)$ is finite.  If $J$ is torsion-free, then as in (1) $J$ has finite $p$-antecedents and we have $J$ is cm-slender and lcH-slender \cite[Theorems A, B(c)]{ConCor}.  Then $\ker(\phi)$ is open and the last case for (2) is complete.

\end{proof}

It should be noted that Baumslag-Solitar groups are themselves known to be n-slender, cm-slender, and lcH-slender \cite[Theorems A, B(i)]{ConCor}.  Finally, we point out that a positive answer to the following question allows one to remove the requirement that the group has no Baumslag-Solitar subgroups.

\begin{question}  If $G$ is a torsion-free one-relator group and $p$ is a prime number greater than the length of the relator of $G$, then does $G$ have finite $p$-antecedents?
\end{question}

\section*{Acknowledgement}

The author thanks the referee for helpful suggestions.

\bibliographystyle{amsplain}

\end{document}